\documentclass[11pt]{article}
\usepackage{latexsym, amsmath, amssymb, amsfonts, amscd}
\usepackage{graphics}
\usepackage[dvips]{epsfig}

\newtheorem{theorem}{Theorem}
\newtheorem{proposition}{Proposition}
\newtheorem{lemma}{Lemma}
\newtheorem{corollary}{Corollary}
\newtheorem{remark}{Remark}
\newtheorem{example}{Example}
\newtheorem{definition}{Definition}

\newtheorem{problem}{Problem}

\newcommand{\complex}{\mathbb C}

\newcommand{\call}{{\cal L}}

\newcommand{\adL}{\mbox{\rm ad}_{\Lambda}}

\newcommand{\abel}{\mbox{\rm{Abel}}}
\newcommand{\kur}{\mbox{\rm{Kur}}}

\def\lcf{\lbrack\! \lbrack}
\def\rcf{\rbrack\! \rbrack}

\def\oell{\overline\ell}
\def\oL{\overline L}
\def\oom{\overline\omega}

\def\ddel{\overline\delta}
\def\dbar{\overline\partial}

\newcommand{\CC}{\mathbb{C}}

\newcommand{\RR}{\mathbb{R}}

\def\ca{{\cal A}}

\def\cg{{\cal G}}
\def\ch{{\cal H}}

\def\cj{{\cal J}}
\def\ck{{\cal K}}

\def\cv{{\cal V}}

\def\dbar{\overline\partial}

\def\om{\omega}
\def\Om{\Omega}

\def\osigma{\overline\sigma}

\def\oomega{\overline\omega}

\def\OLambda{\overline\Lambda}
\def\OG{\overline\Gamma}
\def\otheta{\overline\theta}
\def\osigma{\overline\sigma}

\newcommand{\liea}{\mathfrak{a}}
\newcommand{\lieb}{\mathfrak{b}}

\newcommand{\lieg}{\mathfrak{g}}
\newcommand{\lieh}{\mathfrak{h}}
\newcommand{\liek}{\mathfrak{k}}

\newcommand{\lbra}[2]{\lcf #1, #2 \rcf}

\newcommand{\bproof}{\noindent{\it Proof: }}
\newcommand{\eproof}{\hfill \qed \vspace{0.2in}}
\def\qed{\rule{2.3mm}{2.3mm}}

\begin{document}
\title{\bf Abelian Complex Structures   and Generalizations}
\author{
Yat Sun  Poon\thanks{ Address:
    Department of Mathematics, University of California at Riverside,
    Riverside, CA 92521, U.S.A.. Email: ypoon@ucr.edu.} }
\date{To appear in \it Complex Manifolds \rm}
\maketitle
\begin{abstract} 
After a review on the development of deformation theory of abelian complex structures 
from both the classical and generalized sense, 
we propose the concept of semi-abelian generalized complex structure.  
We present some observations on such structure  and illustrate this new concept with a variety of
examples. 
\end{abstract}

\noindent{AMS Subject Classification: Primary 53D18; Secondary 32G07, 14D15, 17B30, 17B62}

\

\noindent{In memory of my father Yan Ding Poon}


\section{Introduction}

From the onset generalized geometry 
is conceived 
 to encompass complex and symplectic geometry \cite{Marco} \cite{Hitchin-Generalized CY}. 
Its local model theory was settled a few years ago \cite{AB} \cite{Bailey}. 
It adds interests to holomorphic Poisson geometry   \cite{Goto-Poisson} 
\cite{Marco-Pym} \cite{Hitchin-Instanton} \cite{Hitchin-holomorphic Poisson}
  \cite{Pym}. 

While complex structures and symplectic structures are often united in the realm of 
K\"ahler geometry, generalized geometry
enables one to study them with new perspectives, 
especially from the viewpoint of K\"ahler geometry with torsion \cite{Cavalcanti-SKT} \cite{FPS}. 
There is a very rich collection of examples of generalized complex structures
 beyond classical symplectic structures, complex structures, and holomorphic Poisson structures. 
A very prominent class of non-K\"ahlerian manifolds on which complex and symplectic structures coexist is nilmanifolds \cite{BG} 
\cite{Cavalcanti-G} \cite{FPPS}
 \cite{Thurston}. 
 
 This author takes the perspective that generalized geometry, 
especially its deformation theory, is a degree-2 realization of the extended 
deformation theory developed by Kontsevich et al. \cite{B-Kontsevich} 
\cite{CLP} \cite{COP} \cite{Kont} \cite{Merkulov-2003} \cite{Poon-2006}. 
In particular, inspired by the work of Manin and Merkulov
\cite{H-Manin} \cite{Manin-F}
 \cite{Merkulov-2000} \cite{Merkulov-2004},
 this author computed the Frobenius structure on 
 the extended moduli of Kodaira surfaces  \cite{Poon-2006}. 
 One could  restrict an analysis of the Frobenius structure
 from the extended moduli to a generalized moduli in the sense that only degree-2 deformations are allowed. 
 For example, it is proved that the Frobenius structure on the generalized 
moduli of Kodaira manifolds in all dimensions is trivial \cite[Theorem 3]{Poon-2019}. 

Lots of geometric consideration   
and much of the author's work on Frobenius structure and 
holomorphic Poisson deformation  
rely on the underlying complex structures being  abelian
\cite{Andrada-Villacampa}  \cite{Angella-C-Kasuya}
\cite{BS-2000} \cite{BDV} \cite{DF-2000} \cite{Poon-Simanyi-2017} \cite{Poon-Simanyi-2019}.  
This  class of invariant complex structures was first seen in \cite{BDM}, and is often studied 
 along the line of nilpotent complex structures
\cite{CFUG-1996} \cite{CFUG-1997} \cite{Salamon}. 

Some of the  reasons for focusing on abelian complex
structures on nilmanifolds $M$ are due to its accessibility 
in terms of cohomology theory and deformation theory.  
Deformation theory of generalized complex structures is 
based on the concept of Lie bialgebroids \cite{Marco} \cite{LWX} \cite{Mac} \cite{MX}. 
Deformation of a generalized complex structure $\cj$ 
is controlled by its differential Gerstenhaber algebra $DGA(M, \cj)$. 
Gerstenhaber algebra was invented to study deformation of rings after Kodaira's deformation theory
of complex manifolds \cite{Ger}. 
Part of this structure is a Schouten bracket on sections of the exterior algebra of a vector bundle \cite{Mac}. 
When  $\cj$ is a left-invariant abelian complex structure on a
nilmanifold $M=\Delta\backslash G$ modeled on a Lie group $G$,  
$DGA(M, \cj)$ is quasi-isomorphic to its invariant counterpart $DGA(\lieg, \cj)$ \cite{Chen-Fino-P}. 
The invariant $DGA(\lieg, \cj)$ enjoys additional
features as  stated in Proposition \ref{prop: motivation} of the next section. 

We  review the construction of differential Gerstenhaber 
algebra $DGA(M, \cj)$ in Section \ref{sec: gcx}. To motivate our work in Section \ref{sec: semi},
we present  
some  basic properties  when $\cj$ is a classical object such as complex structure, holomorphic Poisson structure, 
and symplectic structure.  

The goal of this paper is to present  the concept of  \it "semi-abelian" \rm generalized complex
structures in terms of  $DGA(\lieg, \cj)$. It is formally stated in
 Definitions \ref{def: key} and \ref{def: key on M}. 
Definition \ref{def: key} 
is made on algebra level. It is designed to capture the features of the invariant part of $DGA(M, \cj)$ as if 
the generalized complex structure $\cj$ is an abelian complex structure. 
It could be extended to analyze generalized complex structures on unimodular algebras \cite{ABD} 
\cite{Andrada-Villacampa}. The definition is made so that all abelian complex structures are semi-abelian (generalized) complex
structures. 

In Section \ref{sec: example}, we provide a collection of examples to illustrate the proposed concept. 
This collection  includes symplectic structure, non-abelian complex structure, and generalized complex
structures of various types. In particular we find the existence of 
semi-abelian (generalized) complex structures that fail to be abelian. There are also nilmanifolds on which there 
are generalized complex structures but none of them could be semi-abelian. 

\section{Abelian Complex Structures}

\subsection{Kodaira surfaces}\label{sec: Kodaira surfaces}
A collection of examples of compact complex surfaces have been inspiring objects in complex manifold theory. 
To name a few, we have complex projective plane, cubic surfaces, K3-surfaces, 
Hopf surfaces, and Kodaira surfaces. 
Kodaira surfaces were discovered as a collection of
elliptic surfaces in classification of compact complex surfaces 
\cite{Kodaira}. They have  trivial 
 canonical bundle. It is also known that the manifold $M$ in question
 is the co-compact quotient of the complex two-dimensional
vector space $\complex^2$ and the co-compact lattice transformation 
leave invariant a complex $(2,0)$-form $\eta$, 
subjected to the conditions
\begin{equation}\label{2 form}
d\eta=0, \qquad \eta\wedge \eta=0, \qquad  \eta\wedge {\overline\eta}>0 
\end{equation}
everywhere on the manifold $M$. As such Kodaira surface is a \it nilmanifold \rm with an invariant complex structure. 

In 1976, Thurston published his well-known paper presenting Kodaira surface as a 
non-K\"ahlerian symplectic manifold \cite{Thurston}
and remarked that there are an abundance of similar examples in higher dimensions. 

Regarding the complex analytic aspect of Kodaira surfaces, 
Borcea examined the moduli of complex structures on Kodaira surfaces by varying 
the co-compact lattices of transformations \cite{Borcea}. 
The starting point of his computation is to recognize that the algebra $\lieg$ on the underlying real vector space for 
$\complex^2$ is solely given by 
\begin{equation}\label{heisenberg}
[X_1, X_2]=2X_4. 
\end{equation}
All other Lie brackets are equal to zero. Therefore, as a Lie algebra it is the direct sum of the trivial one-dimensional 
algebra and the real three-dimensional Heisenberg algebra $\lieh_3$. 
Borcea's computation of the moduli relied on parametrization of the 
complex $2$-form $\eta$ in (\ref{2 form}) 
with respect to the basis of the algebra $\{X_1, X_2, X_3, X_4\}$ and its dual.  

\subsection{Nilpotent complex structures}

On any compact complex surface, the Fr\"olicher spectral sequence associated to the Dolbeault bicomplex degenerates at first step 
\cite{BVP} \cite{Kodaira}. In a cluster of papers \cite{CFG-1987} \cite{CFG-1991} \cite{CFGU-1997} \cite{CFGU-2000}
\cite{CFUG-1996} \cite{CFUG-1997}, Cordero et al. 
studied the degeneracy of Fr\"olicher spectral sequence of invariant complex structures on 
nilmanifolds in all dimensions. 

A compact manifold $M$ is a nilmanifold if it is the quotient of connected 
simply-connected nilpotent Lie group $G$ by a discrete subgroup $\Delta$ so that
$M=\Delta \backslash G$ \cite{Malcev}.  It is a fundamental observation that such 
a Lie group gives rise to a nilmanifold if and only if its
Lie algebra $\lieg$ admits a basis with respect to which the structure constants 
are rational \cite{Malcev}. 
An invariant  complex structure on a nilmanifold $M$ as a right quotient space 
is given by a left-invariant complex structure on the Lie group $G$. 
Equivalently, it is a real linear map
$J: \lieg\to \lieg$  such that $J\circ J=-1$
and its Nijenhuis tensor vanishes. i.e., for all $X, Y$ in $\lieg$, 
\begin{equation}\label{Nijenhuis} 
\lbra{JX}{JY}-\lbra{X}{Y}-J\lbra{JX}{Y}-J\lbra{X}{JY}=0. 
\end{equation}

One of the discoveries by Cordero et al. 
was the 
concept of nilpotent complex structures, see \cite{CFGU-2000} as a preprint reference for \cite{CFG-1987}.   
This subject was further developed by Salamon \cite{Salamon}, 
which has become a key reference for this subject for many publications
including this one. 

Given an invariant complex structure $J$ on $M=\Delta \backslash G$, the space $\lieg^{1,0}$ of $(+i)$-eigenvectors for $J$
in complexified Lie algebra $\lieg_{\CC}$ forms a complex Lie algebra. Denote the $(-i)$-eigenspace by $\lieg^{0,1}$. 
The dual spaces are respectively $\lieg^{*(1,0)}$ and $\lieg^{*(0,1)}$. The complex structure $J$ is nilpotent if there exists
an ordered basis $\{\om^1, \dots, \om^m\}$ for $\lieg^{*(1,0)}$ and constants $A^{j}_{kl}$ and $B^{j}_{kl}$ such that 
for all $1\leq j\leq m$, 
\begin{equation}\label{generic constraint}
d\om^j=\sum_{k, l<j}A^{j}_{kl}\om^k\wedge \oom^l+\sum_{k< l< j}B^{j}_{kl}\om^k\wedge \om^l.
\end{equation}

In the meantime, inspired by a search for geometry with finite holonomy, 
a class of complex structures on nilmanifolds emerged in \cite{BDM}. 
It satisfies the condition that for all $X, Y\in \lieg$, 
\begin{equation}\label{classical abelian}
\lbra{JX}{JY}=\lbra{X}{Y}. 
\end{equation}
This condition is equivalent to require the
 $(+i)$-eigenspace $\lieg^{1,0}$ to form a 
complex abelian algebra although $\lieg$ and $\lieg_{\CC}$ are not necessarily 
abelian. Such complex structures are called \it abelian \rm complex structures. 
In \cite{Salamon} it is proved   that if
a nilpotent algebra $\lieg$ admits an abelian complex structure then there exists an ordered basis 
$\{\om^1, \dots, \om^m\}$ 
for $\lieg^{*(1,0)}$ and constants $A^{j}_{kl}$
 such that 
 for all $1\leq j\leq m$, 
\begin{equation}\label{generic abelian constraint}
d\om^j=\sum_{k, l<j}A^{j}_{kl}\om^k\wedge \oom^l.
\end{equation}
We address this basis for $\lieg^{*(0,1)}$ 
an ascending basis. The observations above lead to the next 
proposition. 

\begin{proposition}\label{prop: motivation} 
Let $J$ be a left-invariant complex structure on a nilmanifold
with Lie algebra $\lieg$. Let  $\lieg^{1,0}$ be the space of invariant 
$(1,0)$-vectors and $\lieg^{*(0,1)}$ the space of invariant $(0,1)$-forms. The following conditions are equivalent. 
\begin{itemize}
\item $\lieg^{1,0}$ is an abelian complex algebra. 
\item $\dbar\oom=0$ for all $\oom\in \lieg^{*(0,1)}$. 
\end{itemize}
$J$ is abelian if and only if one of these conditions is satisfied. 
\end{proposition}

Low-dimension nilpotent algebras are classified  in \cite{Goze}. Given the choice of basis and constraints in 
(\ref{generic constraint}) for nilpotent complex structures and those  
in (\ref{generic abelian constraint}) for abelian complex structures, there is a
classification of the underlying nilpotent algebras to admit such complex structures 
when the real dimension of the algebra is 
at most six  \cite{CFGU-2000} \cite{CFU} \cite{Salamon}. 
A  coarse
dimension count on family of invariant complex structures on each admissible 
algebra was also done in \cite{Salamon}. 
Classification of abelian complex structures in low dimension  is also extended to algebras 
other than the nilpotent ones \cite{ABD}.  

\

\noindent{\bf{Notations.} }
To further present our work, we adopt a convention popularized by \cite{Salamon}. 
Suppose that $e^j, e^k, e^l$ are 1-forms
we use $e^{jkl}$ to represent their exterior product $e^j\wedge e^k \wedge e^l$. 
Likewise when $e_j, e_k$ are 
vectors, the bivector $e_j\wedge e_k$ is represented by $e_{jk}$. 
When $\{e^1, \dots, e^n\}$ is an ordered basis for 
$\lieg^*$, the structure equations for the algebra $\lieg$ in terms of the 
Chevalley-Eilenberg differential is represented
by the n-tuple $(de^1, \dots, de^n)$. 
For example, 
when $de^n=e^{ab}+e^{kl}$  the last entry in  $(de^1, \dots, de^n)$ will be represented by  $ab+kl$. 
\

With the above notations,  the algebra 
of the real 3-dimensional Heisenberg algebra $\lieh_3$ is represented by $(0, 0, 12)$ and 
  the non-trivial algebra with invariant complex structure given in
 (\ref{heisenberg}) is isomorphic to
$(0, 0, 0, 12)=\RR\oplus \lieh_3$. 

\subsection{Examples of abelian complex structures} 
There are only five non-trivial six-dimensional 2-step 
nilpotent algebras admitting abelian complex structures \cite{CFU} \cite{Salamon}. 
Each has a high-dimension generalization. 

\begin{example}\label{ex: Kodaira}  The algebras $(0,0,0,0,0,12)$ and $(0,0,0,0,0,12+34)$.
\end{example}
Consider $(0,0,0,0,0,12)$ as the direct sum $\RR^3\oplus  \lieh_3$. 
The resulting nilmanifold is the product of Kodaira surface with a real two-dimensional torus. 

The four-dimensional algebra $\RR\oplus\lieh_3=(0,0,0,12)$ 
has non-trivial six-dimensional generalization, namely  
$\RR\oplus \lieh_5=(0,0,0,0,0,12+34)$. 
In general, for any natural numbers $m$ and $n$, we have  $\RR^{2m+1}\oplus \lieh_{2n+1}$ 
where $\lieh_{2n+1}$ is the $(2n+1)$-dimensional
Heisenberg algebra with structure equations
$$
(0, \dots, 0,12+\cdots +(2k-1)(2k)+\cdots +(2n-1)(2n))
$$ 
where the first $2n$-entries are all equal to zero. 
An invariant abelian complex structure is  given by $Je_{2k-1}=e_{2k}$ for all $1\leq k\leq n$ and 
$Je_{2n+2m+1}=e_{ 2n+2m+2}.$ The resulting complex manifold corresponding to $\RR\oplus\lieh_{2n+1}$
 is addressed as Kodaira manifolds \cite{Poon-2006}
\cite{Poon-2019}.
We will, however, see the algebra $(0,0,0,0,0, 12+34)$  in a different light in 
Example \ref{ex: 12+34}. 

\begin{example} Three other 2-step six-dimensional nilpotent algebras.
\end{example}
The direct sum of two Heisenberg algebras $(0,0,0,0,12, 34)=\lieh_3\oplus \lieh_3$ 
admits an abelian complex structure $J$ with 
$Je_1=e_2$, $Je_3=e_4$, and $Je_5=e_6$. 
Its high-dimension counterparts $\lieh_{2m+1}\oplus \lieh_{2n+1}$ also have abelian complex structures. 

The algebra $(0, 0, 0, 0, 13+42, 14+23)$ and $(0, 0, 0, 0, 12, 14+23)$ also admit abelian complex structures and 
high-dimension generalization \cite{MPPS-2-step}. 

\begin{example}\label{ex: (0, 0, 0, 12, 14+23, 13+42)}
 A 3-step nilpotent algebra $(0, 0, 0, 12, 14+23, 13+42)$.
\end{example}
After a change of bases $e^4$ to $-e^4$, $e^5$ to $-e^6$ and $e^6$ to $e^5$, we 
present the same algebra with structure equations below. 
\begin{equation}\label{str: (0, 0, 0, 12, 14+23, 13+42)}
de^4=-e^{12}, \qquad de^5=e^{31}+e^{42}, \qquad de^6=e^{41}-e^{32}. 
\end{equation}
Equivalently, 
\[
\lbra{e_1}{e_2}=e_4, \quad
 \lbra{e_3}{e_1}=\lbra{e_4}{e_2}=-e_5, \quad
 \lbra{e_4}{e_1}=-\lbra{e_3}{e_2}=-e_6.
 \]
Define a complex structure $J$ by 
$Je_1=e_2,$ $Je_3=e_4,$ $Je_5=e_6$
 so that $\lieg^{*(1,0)}$ is spanned by $\om^1=e^1+ie^2$, $\om^2=e^3+ie^4$, and $\om^3=e^5+ie^6$. 
As a result of  (\ref{str: (0, 0, 0, 12, 14+23, 13+42)}),
\[
d\om^1=0, \quad d\om^2=\frac12 \om^1\wedge\oom^1,  \quad 
d\om^3=\om^2\wedge \oom^1.
\]  
So, the complex structure $J$ is abelian.

\subsection{Cohomology and deformation} 

Given a nilmanifold $M=\Delta\backslash G$, there is an inclusion of left-invariant differential forms in the space of 
smooth differential forms. 
\begin{equation}
\iota: \wedge^k\lieg^* \longrightarrow C^\infty(M, \wedge^kT^*M).
\end{equation}
Nomizu proved that this inclusion is a quasi-isomorphism in the sense that the inclusion map induces an isomorphism 
at cohomology level  \cite{Nomizu}, 
\[
\iota: H^\bullet (\lieg) \cong H^\bullet_{DR}(M, \RR)
\]
where $H^k(\lieg)$ is the $k$-th cohomology with respect to the complex of the Chevalley-Eilenberg 
differential of $\lieg$. In other words, 
for any $\omega\in \lieg^*$ and $X, Y\in \lieg$, $d\om(X, Y)=-\omega(\lbra{X}{Y}).$ 

Since the attempt by Sakane on similar quasi-isomorphism for Dolbeault cohomology  \cite{Sakane}, 
it has been proved  that 
when $M$ is a nilmanifold with a nilpotent complex structure 
 one obtains  a natural quasi-isomorphism.  
\begin{theorem}{\rm{\cite{CFGU-2000}}}\label{t: Dolbeault}
Suppose that $M=\Delta\backslash G$ is a nilmanifold with a nilpotent complex structure, 
the inclusion of invariant $(p,q)$-forms in the space of sections of $(p,q)$-forms is a quasi-isomorphism. i.e., 
the inclusion map 
\begin{equation}
\iota: \lieg_\CC^{*(p,q)}\longrightarrow C^\infty(M, T^{*(p, q)}M)
\end{equation} 
 induces an isomorphism at cohomology level:
$H^{p,q}_{\dbar}(\lieg_\CC) \cong H^{p,q}_{\dbar}(M). $
\end{theorem}

In addition, Console and Fino initiated a study on the same issue for  all 
 invariant complex structures from the perspective of stability of the desired
  quasi-isomorphisms  \cite{Console} \cite{Console-Fino}. 
It is remarkable that   
the above statement remains an open conjecture when the complex structure is merely invariant. 
For advancement in this direction, 
please see  \cite{Rolle1} \cite{Rolle2} and the latest development in \cite{FRR}  and references therein.

This author and collaborators took on the issue of deformation of invariant complex structures on
 nilmanifolds in \cite{MPPS-2-step}
with a focus on 2-step nilmanifolds with abelian complex structures. The key to enable this investigation was a result 
similar to Theorem \ref{t: Dolbeault}. When $\Theta$ is the sheaf of germs of holomorphic
 tangents for the complex nilmanifold
$M$, the result states as below. 
\begin{theorem}\label{thm: classical}
Suppose that $M=\Delta\backslash G$ is a nilmanifold with an abelian complex structure.
The inclusion map 
\begin{equation}
\iota: \lieg^{1,0}\otimes \lieg^{*(0,k)}\longrightarrow C^\infty(M, T^{1,0}\otimes T^{*(0, k)})
\end{equation} 
induces an isomorphism at cohomology level: 
$H_{\dbar}^k(\lieg^{1,0}) \cong H^k(M, \Theta).$
\end{theorem}

The above theorem was initially proved on 2-step nilmanifolds \cite[Theorem 1]{MPPS-2-step}. It is 
subsequently expanded to include  nilmanifolds with arbitrary number of steps
\cite[Theorem 3.6]{CFP}. Both \cite{CFP} and \cite{MPPS-2-step} 
rely on the fact that the center of the algebra $\lieg$ is invariant of the complex 
structure $J$ when it is abelian. It recreates a (series of) holomorphic fibrations with complex torus as fibers. The proof
of Theorem \ref{thm: classical} becomes an application of Leray spectral sequence formalism.   

The quasi-isomorphism in Theorem \ref{thm: classical}  enables a construction of moduli of complex structures 
on Kodaira manifolds in \cite{GMPP} along the line of Borcea's work in \cite{Borcea}. It also
 enables an analysis on the stability of abelian complex
structures under deformation \cite{CFP}. More broadly, in  \cite{MPPS-2-step} when an abelian complex
structure $J$ is given on a 2-step nilmanifold, the authors considered
 the Kuranishi space ${\kur} (J)$ of the given complex
structure $J$ and the subspace $\abel (J)$ consisting of local deformation 
parameter space of abelian complex structures with 
$J$ in the center. 
Among other results,   they found the following. 
\begin{proposition}{\rm{\cite[Table 1]{MPPS-2-step}}}\label{prop: 6d abelian} 
Consider $M=\Delta\backslash G$ 
a nilmanifold with 2-step nilpotent algebra $\lieg$. 
\begin{itemize}
\item If $\lieg$ is one of  $(0,0,0,0,0,0)$, $(0,0,0,0,0,12)$, and $(0,0,0,0,0, 12+34)$, 
there exists an abelian complex structure $J$ on $M$ such that 
$\dim \abel (J)=\dim \kur(J)$.
\item If $\lieg$ is one of $(0, 0, 0, 0, 12, 34)$, $(0, 0, 0, 0, 13+42, 14+23)$, and $(0, 0, 0, 0, 12, 14+23)$, 
there exists an abelian complex structure
$J$  on the manifold $M$ such that $\dim \abel(J)=\dim \kur(J)-1>0$. 
\end{itemize} 
\end{proposition}

\section{Generalized Complex Structures}\label{sec: gcx}

A generalized complex structure could be conceived as both a tensorial object 
and a spinorial object subjected to multiple conditions. 
Its investigation was initiated by Hitchin \cite{Hitchin-Generalized CY} and developed by Gualtieri \cite{Marco}. 
Given any point $x$ on the manifold $M$, 
when $X, Y$ are in $T_xM$  and 
$\alpha, \beta$ in $T_x^*M$, one considers the pairing between two elements of $T_xM\oplus T_x^*M$
\begin{equation}\label{pairing}
\langle X+\alpha, Y+\beta\rangle =\frac12 (\beta(X)+\alpha (Y)).
\end{equation}
An almost generalized complex structure is a bundle map 
\begin{equation}
\cj: TM\oplus T^*M \to TM\oplus T^*M
\end{equation} 
such that $\cj\circ \cj=-1$ and 
$
 \langle \cj(X+\alpha), \cj(Y+\beta)\rangle= \langle X+\alpha, Y+\beta\rangle
$
for all $X+\alpha, Y+\beta$ in  $T_xM\oplus T_x^*M$ for each point $x$ in the manifold $M$. 
It is known that it could be expressed as a map of bundle of direct sum $TM\oplus T^*M$, 
\begin{equation}\label{matrix rep}
\cj=
\left(
\begin{array}{cc}
J & \Pi \\
B & -J^*
\end{array}
\right)
\end{equation}
where $B$ is a two-form and $\Pi$ is a bivector \cite{Marco}. 
Equivalently, the bundle of $(+i)$-eigenspace $L$ is maximally isotropic subbundle of the complexification of
$TM\oplus T^*M$. Let $\oL$ represent the complex conjugate bundle of $L$, it is the bundle of $(-i)$-eigenspace. 
By virtual of $L$ being maximally isotropic,  $\oL$ is also maximally isotropic and $L\oplus \oL =(TM\oplus T^*M)_\CC$. Since 
$L$ and $\oL$ are maximally isotropic and the bilinear form $\langle -,-\rangle$ is non-degenerate, it yields naturally isomorphisms 
by the pairing  ${\langle -,-\rangle}$ as given in (\ref{pairing}): 
\begin{equation}
L\cong {\oL}^*, \qquad \oL \cong L^*.
\end{equation}
In subsequent discussion, we often utilize these isomorphisms. 

The Courant bracket defined on the space of sections $C^\infty(M, TM\oplus T^*M)$ is 
\begin{equation}
\lbra{X+\alpha}{Y+\beta}=\lbra{X}{Y}+\call_X\beta-\call_Y\alpha-\frac12 d(\beta(X)-\alpha (Y)),
\end{equation}
where $\call_X\beta$ is the Lie derivative of the one-form $\beta$ along the vector field $X$ and 
$\lbra{X}{Y}$ is the Lie bracket of vector fields. 

An almost generalized complex structure is integrable if the space of smooth 
sections of the bundle $L$, denoted by $C^\infty(M, L)$, is closed with respect to the Courant bracket. 
Although the Courant bracket on $C^\infty(M, TM\oplus T^*M)$ does not satisfy the Jacobi identity, 
it is observed in \cite[Proposition 3.27]{Marco} that its restriction on $C^\infty(M, L)$ satisfies the 
Jacobi identity if and only if $C^\infty(M, L)$ is closed under the Courant bracket. The fact that 
$L$ is maximally isotropic plays a key role in this observation. 

\subsection{DGA and cohomology}
Under the assumption that a generalized complex structure $\cj$ is integrable, one  treats $L$ as a 
Lie algebroid \cite{Marco} 
\cite{LWX} and extends the restriction of the Courant bracket to the space of section of exterior algebras
$C^\infty(M, \wedge^\bullet L)$ to obtain what is known as a Gerstenhaber algebra, also as
a Schouten algebra \cite{Ger} \cite[Definition 7.5.1]{Mac}. We will address the restriction of the Courant 
bracket on $C^\infty(M, \wedge^\bullet L)$  as the Schouten bracket. 

Taking complex conjugations, the space $C^\infty(M, \wedge^\bullet \oL)$ also inherits a Schouten bracket.  
Since $\oL$ is a Lie algebroid, it has an associated differential on its dual ${\oL}^*\cong L$:
\begin{equation}
\ddel: C^\infty(M, L) \longrightarrow C^\infty(M, \wedge^2L).
\end{equation}
It is  observed in \cite{Marco} that the pair $(L, \oL)$ forms a Lie bialgebroid in the sense of 
\cite{MX}. 
The operator $\ddel$ extends to an even exterior differential operator and 
an odd differential Lie superalgebra operator \cite{Manin} \cite{Poon-2006}. It means that
if $a\in C^\infty(M, \wedge^{|a|}L)$ and $b\in C^\infty(M, \wedge^{|b|}L)$, then 
\begin{equation}\label{compatible differential} 
\ddel (a\wedge b)=(\ddel a)\wedge b+(-1)^{|a|} a\wedge(\ddel b), 
\qquad
\ddel\lbra{a}{b}=\lbra{\ddel a}{b}-(-1)^{|a|}\lbra{a}{\ddel b}. 
\end{equation} 
In summary, each generalized complex structure is associated to  
\begin{equation}\label{symbolic DGA}
DGA(M, \cj)=\left( C^\infty(M, \wedge^\bullet \oL), \wedge, \lbra{\bullet}{\bullet}, \ddel \right),
\end{equation}
 a differential Gerstenhaber algebra \cite{Ger} \cite{Poon-2006}. 

It is further observed in \cite{Marco} that the operator $\ddel$ is elliptic and the integrability of 
$\cj$ is equivalent to $\ddel\circ \ddel=0$. Therefore $DGA(M, \cj)$ has an associated cohomology theory
\[
H^k(M, \cj)=\frac{\ker \ddel: C^\infty(M, \wedge^kL) \longrightarrow C^\infty(M, \wedge^{k+1}L)}
{\mbox{\rm{image} } \ddel: C^\infty(M, \wedge^{k-1}L) \longrightarrow C^\infty(M, \wedge^{k}L)}.
\]
As a consequence of (\ref{compatible differential}),  the exterior product 
and the Schouten bracket descend to the cohomology space
so that we obtain a Gerstenhaber algebra:
\[
\left( H^\bullet (M, \cj), \wedge, \lbra{\bullet}{\bullet}\right).
\]

From the view point of extended deformation \cite{Kont} \cite{Merkulov-2003} \cite{Poon-2006}, one treats an 
element $\Gamma$ in $C^\infty(M, \wedge^\bullet L)$ as an integrable \it extended \rm deformation 
if it satisfies the Maurer-Cartan equation 
\begin{equation}\label{MC}
 \ddel \Gamma+\frac12\lbra{\Gamma}{\Gamma}=0.
\end{equation}
We obtain a generalized deformation 
when $\Gamma$ is a degree-2 section, i.e., an element in $C^\infty(M, \wedge^2 L)$.

\subsection{Symplectic manifolds} 
Suppose that $M$ is a smooth manifold and $\Omega$ is a symplectic form, 
we  treat a contraction of $\Omega$ with a tangent vector
as a bundle map $\Omega: TM\to T^*M$. Since $\Omega$ is non-degenerate, its inverse $\Omega^{-1}$ is 
a well-defined 
 bivector. In the matrix representation of a generalized complex structure $\cj$, 
by choosing $J=0$, $B=\Omega$, and $\Pi=\Omega^{-1}$, one obtains a generalized complex structure.  
Pointwisely, 
\begin{equation}
L=\{X-i\Om (X): X\in TM\}, \qquad \oL=\{X+i\Om (X): X\in TM\}.
\end{equation}
The integrability of $\cj$ as a generalized complex structure amounts to the identity
\begin{equation}\label{symp 1}
\lbra{X-i\Om (X)}{Y-i\Om (Y)}=\lbra{X}{Y}-i\Om(\lbra{X}{Y})
\end{equation}
for all $X, Y$ in $C^\infty(M, TM)$. It is equivalent to $d\Omega=0$.
Let $X, Y, Z$ be any smooth vector fields on $M$. 
\begin{eqnarray}
&& \ddel (Z-i\Omega(Z))({X+i\Omega(X)},{Y+i\Omega(Y)}) \nonumber\\
&=&-\langle Z-i\Omega(Z), \lbra{X+i\Omega(X)}{Y+i\Omega(Y)}\rangle 
=-\langle Z-i\Omega(Z), \lbra{X}{Y}+i\Omega(\lbra{X}{Y})  \rangle 
\nonumber 
\\
&=&i (\Omega(Z))\lbra{X}{Y}=-i (d(\Omega(Z)))({X},{Y}). \label{symp 2}
\end{eqnarray}
Consider a bundle map 
\begin{equation}\label{bundle map}
\varphi: L \hookrightarrow (TM\oplus T^*M)_\CC \longrightarrow T^*M_\CC
\end{equation}
obtained by a composition of inclusion and projection. 
$\varphi$ is a bundle isomorphism because $\Omega$ is non-degenerate. 
It extends naturally to a bundle map of exterior powers
so that for all $a, b\in C^\infty(M,  L)$,  $\varphi(a)\wedge \varphi(b)=\varphi(a\wedge b)$. 
By (\ref{symp 2})
$\varphi( \ddel a) =d (\varphi(a))$. Moreover, for any $\mu, \nu$ in $C^\infty(M,  T^*M_\CC)$ define
\begin{equation}
\lbra{\mu}{\nu}_\Omega:=\lbra{\Omega^{-1}\mu}{\Omega^{-1}\nu}
\end{equation}
where the bracket on the right hand side is the usual Lie bracket of vector fields. By  
(\ref{symp 1}), $\lbra{a}{b}=\lbra{\varphi(a)}{\varphi(b)}_\Omega.$ 
In particular, the differential Gerstenhaber algebra 
$DGA(M, \Omega)$ defined on $C^\infty(M, \wedge^\bullet L)$ with Courant bracket
$\lbra{\bullet}{\bullet}$ 
 is isomorphic to 
the differential Gerstenhaber algebra defined on the space of differential forms
\begin{equation}
\left( C^\infty(M, \wedge^\bullet T^*M_\CC), \wedge, \lbra{\bullet}{\bullet}_\Omega, d \right)
\end{equation}
where the differential  $d$  is  the deRham differential. The resulting cohomology is the 
complexified deRham cohomology
$(H^\bullet_{DR}(M, \CC), \wedge, \lbra{\bullet}{\bullet}_\Omega). $

\begin{proposition}\label{proposition: symp iso} When $\Omega$ is a symplectic form, the bundle map {\rm (\ref{bundle map})} defines an 
 isomorphism of differential Gerstenhaber algebra. 
\begin{equation}\label{symp iso}
\varphi: DGA(M, \Omega) \cong
\left( C^\infty(M, \wedge^\bullet T^*M_\CC), \wedge, \lbra{\bullet}{\bullet}_\Omega, d\right)
\end{equation}
\end{proposition}

The isomorphism above demonstrates that the differential Gerstenhaber algebra of a symplectic structure
 as a generalized complex structure is consistent with  the ones in 
  extended deformation theory in \cite{Kont} \cite{Merkulov-2000} 
\cite{Merkulov-2003} \cite{Poon-2006} and weak mirror symmetry 
\cite{CLP} \cite{COP} \cite{Merkulov-2003}.

\subsection{Classical complex structures}
The matrix representation of $\cj$ in (\ref{matrix rep}) for a classical complex structure 
is given by $B=0$ and $\Pi=0$. 
The bundles of $(+i)$ and $(-i)$ eigenvectors with respect to $\cj$ are respectively 
\begin{equation}\label{L for complex}
L=T^{1,0}\oplus T^{*(0,1)}, \qquad 
\oL=T^{0,1}\oplus T^{*(1,0)}
\end{equation}
where $T^{1,0}$ and $T^{0,1}$ are  bundles of (1,0)-vectors and (0,1)-vectors  with respect to $J$. 
$T^{*(1,0)}$ and $T^{*(0,1)}$ are their respective duals. 

\begin{proposition} \label{prop: classical}
Let $\cj$ be a classical complex structure and $\ddel$ the Lie algebroid differential on $C^\infty(M, T^{1,0}\oplus T^{*(0,1)})$. 
For any (1,0)-vector field $U$, $\ddel U=\frac12 \dbar U$.
For any (0,1)-form $\oom$, $\ddel \oom=\frac12 \dbar\oom$. 
\end{proposition}
\bproof Given $U, V, W\in C^\infty(M, T^{1,0})$ and $\om, \mu, \nu\in C^\infty(M, T^{*(0,1)})$, we  test 
$\ddel U$ and $\ddel \oom$  against the pairs $\{{\overline V}, {\overline W}\}$, $\{ {\overline V},{\nu} \}$, and 
$\{ {\mu},{\nu} \}$ respectively. As $\lbra{\mu}{\nu}=0$, we have only two pairs to work on. 
\begin{eqnarray*}
&\ddel U( {\overline V}, {\overline W}) =-\langle U, \lbra{\overline V}{\overline W} \rangle=0, &\\
&\ddel U({\overline V},{\nu}) =-\langle U, \lbra{\overline V}{\nu}\rangle
=-\frac12 \lbra{\overline V}{\nu}(U)= -\frac12 d\nu({\overline V}, U)=\frac12 \nu (\lbra{\overline V}{U}). &
\end{eqnarray*}
Since $\nu$ is a (1,0)-form, $\nu (\lbra{\overline V}{U})=\nu (\lbra{\overline V}{U}^{1,0})=\nu(\dbar_{\overline V}U)$ 
\cite{Gau} \cite{Poon-2006}. We obtain $\ddel U=\frac12 \dbar U$. 
Similarly $\ddel \oom({\overline V},{\nu})=0$ and 
$$
\ddel \oom( {\overline V}, {\overline W}) =-\langle \oom, \lbra{\overline V}{\overline W} \rangle
=\frac12 d\oom ({\overline V},{\overline W} )=\frac12 \dbar\oom ({\overline V},{\overline W} ).
 $$
 It follows that $\ddel \oom=\frac12 \dbar\oom$. 
\eproof

The observation above implies  that a (1,0)-vector field $U$ is holomorphic if and only
 if $\ddel U=0$. Equivalently, 
the (1,0)-component of $\lbra{\overline V}{U}$ vanishes for any (0,1)-vector field $\overline V$. 
For instance, in terms of the ascending basis of a nilpotent 
complex structure as given in 
(\ref{generic abelian constraint}) and its dual basis $\{T_1, \dots, T_m\}$, $\dbar T_m=0$. 
For further reference, we have
\begin{corollary}\label{holomorphic Poisson} 
Let $\Lambda$ be a (2,0)-vector field. It is holomorphic $\dbar\Lambda=0$ 
if and only if for any (0,1)-vector field $\overline V$, the (2,0)-component of 
$\lbra{\overline V}{\Lambda}$ vanishes. 
\end{corollary}

As a result, it is now also obvious that in terms of the ascending basis and its dual $\{T_1, \dots, T_m\}$, 
$\dbar({T_{m-1}\wedge T_m})=0$ and $\lbra{T_{m-1}\wedge T_m}{T_{m-1}\wedge T_m}=0$, a fact observed 
in \cite[Theorem 5.1]{Cavalcanti-G}.

\subsection{Holomorphic Poisson structures}
Holomorphic Poisson structure  is built upon a complex structure $J$. 
In a matrix representation for $\cj$, $J$ is a complex
structure, $B=0$, and $\Pi\neq 0$. In the presence of  $J$, it is convenient to 
put $\Pi$ in complex terms. A contraction with $\Pi$ defines a linear map 
$
\Pi: T^*M \to TM.
$
The complex structure $J$ maps from $TM$ to $TM$ while $J^*$ maps from $T^*M$ to $T^*M$. Define
$
\Upsilon: T^*M \to TM$  by $\Upsilon=-\Pi\circ J^*.$
It follows that $\Pi=\Upsilon\circ J^*$ because $J\circ J=-1$.  $\Pi+i\Upsilon$ is a type-(2,0) bivector 
because its contraction with any (0,1)-form is equal to zero. Define
\begin{equation}
\Lambda=\frac{i}2(\Pi+i\Upsilon).
\end{equation} 
The $(+i)$-eigenbundle $L$ with respect to $\cj$ has a real representation in terms of  $\Pi$ 
and a complex representation in terms of $\Lambda$. 
\begin{eqnarray*}
&L=\{ X-iJX, \alpha+iJ^*\alpha-i\Pi(\alpha): X\in TM, \alpha\in T^*M\}, &\\
&L=\{U, \oom+\OLambda\oom: U\in T^{1,0}, \oom\in T^{*(0,1)}\}.& \label{L for Poisson}
\end{eqnarray*}
 It is most convenient to express the integrability of $L$ in terms of the second representation because its 
 integrability is equivalent to $J$ being an integrable classical complex structure and $\Lambda$ is a
 holomorphic Poisson structure, i.e., $\dbar\Lambda=0$ and $\lbra{\Lambda}{\Lambda}=0$. 

\begin{lemma}\label{lemma: ddel u} For any (1,0)-vector field $U$ on a holomorphic Poisson manifold, 
 $\ddel U=0$ if and only if  $\dbar U=0$ and $\call_U\Lambda=0$, i.e., $U$ is an infinitesimal 
holomorphic Poisson transformation.  
\end{lemma}
\bproof
Note that for any (0,1)-vector fields $\overline V$ and $\overline W$, 
$(\ddel U)(\overline{V}, \overline{W})=0$
because $TM_{\CC}$ is isotropic. When $\om$ is a (1,0)-form, 
\begin{eqnarray*}
&& (\ddel U)(\overline{V},{\om+\Lambda\om})=-\langle U, \lbra{\overline{V}}{\om+\Lambda\om}\rangle 
= -\langle U, \lbra{\overline{V}}{\Lambda\om}+\lbra{\overline{V}}{\om}\rangle\\
&=& -\langle U, \lbra{\overline{V}}{\om}\rangle
=-\frac12 (\iota_{\overline{V}}d\om)(U)=-\frac12 d\om(\overline{V}, U)=\frac12 \om(\lbra{\overline{V}}{U}). 
\end{eqnarray*}
The last is equal to zero for all $\om$ and all $\overline{V}$ if and only if $\dbar U=0$ \cite{Gau}. Next
for any (1,0)-forms $\mu, \nu$, 
\begin{eqnarray*} 
&&(\ddel U)({\mu+\Lambda\mu},{\nu+\Lambda\nu})=-\langle U, \lbra{\mu+\Lambda\mu}{\nu+ \Lambda\nu}\rangle\\
&=&
\langle U, \lbra{\Lambda\mu}{\nu} -\lbra{\Lambda\nu}{\mu} \rangle
=-\frac12\left( d\nu(\Lambda\mu, U)-d\mu(\Lambda\nu, U)  \right)\\
&=&\frac12\left( \nu(\lbra{\Lambda\mu}{U})-\mu(\lbra{\Lambda\nu}{U})    \right)=-\frac12 (\call_U\Lambda)(\mu, \nu)
=\frac12\lbra{\Lambda}{U}(\mu, \nu).
\end{eqnarray*}
It is equal to zero for all $\mu, \nu$ if and only if $\call_U\Lambda=0$. 
\eproof
\begin{lemma}\label{lemma: ddel oom}  For any (0,1)-form $\oom$,     
 $\ddel ( \oom+{\overline\Lambda}\oomega)=0$ if and only if 
   $\lbra{\Lambda}{\oom}=0$ and $ \dbar\oom=0$. 
 \end{lemma}
 \bproof
 For any (0,1)-vector fields ${\overline V}$ and ${\overline W}$, 
\begin{eqnarray*}
&&\ddel ( \oom+{\overline\Lambda}\oomega)(\overline{V}, \overline{W})
=-\langle \oom+{\overline\Lambda}\oomega, \lbra{\overline{V}}{\overline{W}}\rangle
\\
&=&
-\frac12 \oom(\lbra{\overline{V}}{\overline{W}})=\frac12 (d\oom)({\overline{V}},{\overline{W}})
=\frac12 (\dbar\oom)({\overline{V}},{\overline{W}}).
\end{eqnarray*}
Let $\nu$ be any (1,0)-form, 
\begin{eqnarray*}
&&\ddel ( \oom+{\overline\Lambda}\oomega)(\overline{V},{\nu+\Lambda\nu})\\
&=&-\langle \oom+{\overline\Lambda}\oomega, \lbra{\overline{V}}{\nu+\Lambda\nu}\rangle
=-\langle {\overline\Lambda}\oomega+\oom, \lbra{\overline{V}}{\Lambda\nu}+\lbra{\overline{V}}{\nu}\rangle\\
&=&-\frac12\left( 
\oom( \lbra{\overline{V}}{\Lambda\nu})+ \lbra{\overline{V}}{\nu}(\OLambda\oom)
\right)
=-\frac12\left( 
\oom( \lbra{\overline{V}}{\Lambda\nu})+ (\iota_{\overline{V}}{d\nu})(\OLambda\oom)
\right).
\end{eqnarray*}
As the complex structure is integrable, the type-(0,2) component of $d\nu$ vanishes. Therefore, 
\begin{equation}
\ddel ( \oom+{\overline\Lambda}\oomega)(\overline{V},{\nu+\Lambda\nu})=
-\frac12 \oom( \lbra{\overline{V}}{\Lambda\nu})=\frac12  \lbra{\Lambda}{\oom}(\nu, {\overline V}).
\end{equation}
Next, for any (1,0)-forms $\mu, \nu$, 
\begin{eqnarray*}
&&\ddel ( \oom+{\overline\Lambda}\oomega)({\mu+\Lambda\mu}, {\nu+\Lambda\nu})
=-\langle  \oom+{\overline\Lambda}\oomega, \lbra{\mu+\Lambda\mu}{\nu+\Lambda\nu} \rangle\\
&=&-\frac12\left(
\oom( \lbra{\Lambda\mu}{\Lambda\nu} )+\lbra{\Lambda\mu}{\nu}({\OLambda}\oom)
-\lbra{\Lambda\nu}{\mu}({\OLambda}\oom)
\right)\\
&=&-\frac12\left(\lbra{\Lambda\mu}{\nu}({\OLambda}\oom)
-\lbra{\Lambda\nu}{\mu}({\OLambda}\oom)
\right)
=-\frac12\left(
\OLambda(\oom, \lbra{\Lambda\mu}{\nu})-\OLambda(\oom, \lbra{\Lambda\nu}{\mu})
\right)\\
&=&\frac12 \lbra{\Lambda}{\OLambda\oom}(\mu, \nu).
\end{eqnarray*} 
By Corollary \ref{holomorphic Poisson}, when $\Lambda$ is holomorphic
the $(2,0)$-component of $\lbra{\Lambda}{\OLambda\oom}$  vanishes for any (0,1)-form $\oom$. 
Therefore, the only obstructions for 
$ \oom+{\overline\Lambda}\oomega$ to be $\ddel$-closed is when $\dbar\oom=0$ 
and $\lbra{\Lambda}{\oom}=0$ as stated.
\eproof

\subsection{Different types of generalized complex structures}\label{sec: type}
There are many examples of generalized complex structures different from the symplectic 
or complex types.
For example, one could have a symplectic bundle
over complex manifold \cite{Angella-C-Kasuya}. 
On the total space of such a fiber bundle, 
there exists a closed 2-form such that its restriction to each fiber is a symplectic form. The base manifold
is a complex manifold. Below we apply a well-known method to construct a 
generalized complex structure on a manifold that is neither symplectic type or complex type.  
Let $U(1)$ represent the one-dimensional unitary group. 

\begin{proposition}\label{prop: toric} Suppose that $M$ is the total space of a principal $U(1)\times U(1)$-bundle 
over a complex n-dimensional manifold $X$ . If the  curvature of a connection on this
bundle is represented by type-(1,1) forms on $X$,
 then $M$ admits a generalized complex structure of type-n on the manifold $M$. 
 \end{proposition} 
 \bproof
We prove this theorem by mimicking a well-known construction of integrable complex structure on toric bundles. 
See e.g., \cite{GGP}. 

Denote the principal bundle projection by $\pi:M\to X$. 
Let $\theta=(\theta^1, \theta^2): M\to \RR\oplus \RR$ be the connection 1-form.
The curvature form of this connection is 
$(d\theta^1, d\theta^2)$. By assumption, there are type-(1,1) forms $\vartheta^1, \vartheta^2$ on the manifold $X$ such that 
$d\theta^1=\pi^*\vartheta^1$ and $d\theta^2=\pi^*\vartheta^2$. 
Let $\cv$ be the bundle of vertical vector fields generated by 
the principal action of $U(1)\times U(1)$. It yields an exact sequence
of vector bundles  on the manifold $M$:
$$
0\to \cv\hookrightarrow TM \stackrel{d\pi}{\longrightarrow} \pi^*TX\to 0.
$$
Let  ${\cal{H}}=\ker\theta=\ker\theta^1 \cap\ker\theta^2$ be 
the horizontal space of the connection $\theta$. It yields a splitting: $TM=\ch\oplus \cv$. The projection 
$d\pi: TM\to \pi^*TX$ yields an isomorphism  $d\pi: {\cal{H}}\to \pi^*TX$ via horizontal lift. For any tangent
vector  $e$ on the manifold $X$, denote its horizontal lift to $M$ by $e^h$ so that $d\pi(e^h)=e$.

Consider $B=\theta^1\wedge \theta^2$ as a 2-form on the 
manifold $M$. Let $(v_1, v_2)$ be the fundamental vector field on
$M$ generated by the principal $U(1)\times U(1)$-action. In particular, they  are vertical vector fields
trivializing the vertical  bundle $\cv$ over the manifold $M$. Fiberwise, they form a dual to the vertical 
1-forms $(\theta_1, \theta_2)$. 
Let $\Pi=v_1\wedge v_2$ be
a bivector field on the manifold $M$. 
We now obtain two of the three components of the matrix representation  (\ref{matrix rep}) of 
an (almost) generalized complex structure 
$\cj$ on $M$. The last component $J$ acts on $\ch$ and it is equal to zero on $\cv$. 

   Define $J$ on $\cal{H}$ by the horizontal lift of
the action of $J$ on the base manifold $X$. i.e.,  for each $e^h$ local section of $\cal{H}$, 
\begin{equation}
d\pi (J e^h)=J d\pi (e^h).
\end{equation}
Therefore, we obtain an \it almost \rm generalized complex structure $\cj$ with $(+i)$-eigenbundles given by
\[
L ={\cal{H}}^{1,0}\oplus\pi^*T^{*(0,1)}X\oplus \{v-i B(v),  \rho-i \Pi (\rho): v \in \cv, \rho\in \cv^*\}.
\]
It is apparent that $L$ is isotropic with respect to the non-degenerate pairing (\ref{pairing}). The space 
$ \{v-i B(v),  \rho-i \Pi (\rho): v \in \cv, \rho\in \cv^*\}$ is  trivialized by 
\begin{equation}
 v_1-i B(v_1)=v_1-i\theta^2, \quad  v_2-i B(v_2)=v_2+i\theta^1.
 \end{equation}

If $u, u_1, u_2$ are (1,0)-vector fields on the base manifold $X$,  
$u_1^h$ and $u_2^h$ are $(+i)$-eigenvectors with respect to $\cj$ on 
the manifold $M$. On  $M$, 
$$
\lbra{u_1^h}{u_2^h}=(d\theta^1, d\theta^2)(u_1^h, u_2^h)=(\pi^*\vartheta^1, \pi^*\vartheta^2)(u_1^h, u_2^h)
=(\vartheta^1(u_1, u_2), \vartheta^2(u_1, u_2))=0
$$ 
because $\vartheta^1$ and $\vartheta^2$ are both type-(1,1) forms. 

On the other hand, a vertical $(+i)$-eigenvector is given by 
$v_1-iB(v_1)=v_1-i\theta^2$. Since all horizontal distributions are invariant
of the principal action and $u^h$ is in $\ker\theta^1\cap\ker\theta^2$, 
\begin{eqnarray*}
&& \lbra{u^h}{v_1-iB(v_1)}= \lbra{u^h}{v_1-i\theta^2}=-i\lbra{u^h}{\theta^2} \\
&=& -i( \call_{u^h}{\theta^2}-\frac12 d(\theta^2 (u^h) ) )
=-i\iota_{u^h}d\theta^2=-i\iota_{u^h}\pi^*\vartheta^2=-i \pi^*(\iota_u\vartheta^2).
\end{eqnarray*}
Since $\vartheta^2$ is type-(1,1) and $u$ is a type-(1,0), $\iota_u\vartheta^2$ is a type-(0,1) form
on the manifold $X$. It follows that last term is a section of $\pi^*T^{*(0,1)}X\subset L$.
If $\oom$ is a (0,1)-form on the base complex manifold $X$, 
\[
\lbra{\pi^*\oom}{v_1-iB(v_1)}= \lbra{\pi^*\oom}{v_1-i\theta^2}=\lbra{\pi^*\oom}{v_1} 
=-\call_{v_1}(\pi^*\oom)+\frac12 d((\pi^*\oom)(v_1)).
\]
Since the pull-back form $\pi^*\oom$ is invariant of vertical vector field, the above is equal to zero. 

A similar computation works for $\lbra{u^h}{v_2-iB(v_2)}=\lbra{u^h}{v_2+i\theta^1}$. Therefore, 
the space of sections for $L$ is closed with respect to the Courant bracket. 
\eproof

There are lots examples of non-trivial toric bundles as described by Proposition \ref{prop: toric} with very interesting 
complex geometry; see e.g., \cite{GGP}. 
In this note, we illustrate the above construction with Example \ref{ex: 12+34}.

\subsection{Deformation} 
Recall that the construction for the differential
Gerstenhaber algebra $DGA(M, \cj)$ for a generalized complex structure $\cj$ on a smooth manifold $M$ is grounded on the 
Lie bialgebroid $(L, \oL)$ of $(+i)$ and $(-i)$-eigenbundles with respect to bundle map $\cj$. The differential $\ddel$ on $L$ is 
a reflection of the structure of the Schouten bracket on sections of $\wedge^\bullet\oL$ and $L$ is treated as the dual bundle of 
$\oL$ via the natural non-degenerate pairing. Given  $\Gamma\in C^\infty(M, \wedge^2L)$, we treat it as a map
$\Gamma: \oL=L^*\to L$ and its conjugation as another map $\OG: L=\oL^* \to \oL $.
Define
\begin{equation}
L_{\OG}=\{\sigma+\OG(\sigma): \sigma\in L\}, \qquad \oL_\Gamma=\{\osigma+\Gamma(\osigma):    \osigma\in \oL\}.
\end{equation}
The pair forms a new Lie bialgebroid, or equivalently new generalized complex structure in our context if and only if 
$\Gamma$ satisfies the Maurer-Cartan equation (\ref{MC}) and the pair stays maximally isotropic \cite{Marco} \cite{LWX}.  
Furthermore, the direct sum $ L_{\OG}\oplus \oL_\Gamma$ has alternative representation when $\Gamma$ is sufficiently 
close to zero. 
 \begin{equation}
 L_{\OG}\oplus \oL_\Gamma=(TM\oplus T^*M)_\CC=L\oplus \oL_\Gamma.
 \end{equation}
 From the right hand side, the space of sections for $\wedge^\bullet L$ continues to inherit the
 Courant bracket on $(TM\oplus T^*M)_\CC$. However, the natural non-degenerate pairing 
 identifies $L$ to the dual of $\oL_\Gamma$. Therefore
 $
 \wedge^\bullet L=\wedge^\bullet (\oL_\Gamma)^*, 
 $
 and it
 inherits a differential $\ddel_\Gamma$ determined by the restriction of the 
 Courant bracket  on  $\wedge^\bullet \oL_\Gamma$. 
 
 \begin{theorem}{\rm\cite{LWX}}\label{thm: lwx}
The differential $\ddel$ on $L$ as a dual to $\oL$ and the differential 
 $\ddel_\Gamma$ on $L$ as a dual to $ \oL_\Gamma$ are related by
\begin{equation}
\ddel_\Gamma a=\ddel a+\lbra{\Gamma}{a}
\end{equation}
for all section $a$ in $C^\infty(M, \wedge^\bullet L)$. In particular,  $DGA(M, \cj_\Gamma)$ after deformation 
by $\Gamma$ is
isomorphic to $(C^\infty(M, \wedge^\bullet L), \wedge, \lbra{-}{-}, \ddel_\Gamma)$.
\end{theorem}

From now on, we denote $\lbra{\Gamma}{a}$ by $ad_\Gamma(a)$. In the notations above, the Maurer-Cartan 
equation (\ref{MC}) is translated to $\ddel_\Gamma\circ\ddel_\Gamma=0$ \cite{LWX}. The next corollary is trivial. 

\begin{corollary}\label{thm: trivial} Let $\cj$ be a generalized complex structure  on a 
manifold  $M$  with differential 
Gerstenhaber algebra 
$DGA(M, \cj)$.  
If $\Gamma$ is an element in $C^\infty(M, \wedge^2L)$ such that 
$\ddel \Gamma=0$  and $ad_\Gamma=0.$
Then $DGA(M, \cj)=DGA(M, \cj_\Gamma)$.
\end{corollary}

When a generalized complex structure $\cj$ is given by a holomorphic Poisson
structure $(J, \Lambda)$, from the viewpoint of  Theorem \ref{thm: lwx} 
we treat $\Lambda$ as a deformation of the complex structure $J$. 
In terms of a matrix representation with components $J$ and $\Pi$, we have 
\begin{equation}
\cj_t=
\left(
\begin{array}{cc}
J & t\Pi \\
0 & -J^*
\end{array}
\right)
\end{equation}
where $t$ is the deformation parameter. When $t=0$, we reach the underlying complex structure. 
When $t=1$, we have our holomorphic Poisson structure. In this
perspective, the algebroid differential $\ddel$ with respect to the generalized
complex structure associated to $(J, \Lambda)$ is identified to $\dbar+\adL$ acting on the 
Lie algebroid of the complex structure $J$, and Lemma \ref{lemma: ddel u} 
and Lemma \ref{lemma: ddel oom} 
follow easily. 

\begin{example} A holomorphic Poisson structure associated to an abelian complex structure on the
 algebra $(0,0,0, 12, 14+23, 13+42)$.
\end{example}
We revisit Example \ref{ex: (0, 0, 0, 12, 14+23, 13+42)}. Recall that
$Je_1=e_2$, $Je_3=e_4$, $Je_5=e_6.$ 
 Let $T_1=\frac12(e_1-ie_2)$, 
$T_2=\frac12(e_3-ie_4)$, and $T_3=\frac12(e_5-ie_6)$. Then $\{T_1, T_2, T_3\}$ spans $\lieg^{1,0}$. Now we could work out 
the structure equations for $DGA(\lieg, J)$. 
\begin{equation}\label{poisson example structure}
d\oom^1=0, \qquad d\oom^2=-\frac12\om^1\wedge\oom^1, \qquad
d\oom^3=-\om^1\wedge\oom^2.
\end{equation}
 In particular, the non-zero Courant brackets on $\lieg_\CC$ are
\begin{equation}
\lbra{T_1}{ {\overline T}_1}=-\frac12 T_2+\frac12 {\overline T}_2, \quad 
\lbra{T_1}{ {\overline T}_2}={\overline T}_3, \quad 
\lbra{T_2}{ {\overline T}_1}=-T_3. 
\end{equation}
Since $\lbra{T_j}{\oom^k}=\iota_{T_j}d\oom^k$, the non-zero brackets between 
elements in $\lieg^{1,0}$ and elements in 
$\lieg^{*(0,1)}$ are
\begin{equation}
\lbra{T_1}{\oom^2}=-\frac12 \oom^1, \quad \lbra{T_1}{\oom^3}=-\oom^2. 
\end{equation}
Finally, as for any $T\in \lieg^{1,0}$, $\dbar T=\sum_{k=1}^3\lbra{ {\overline T}_k}{T}^{1,0}\wedge \oom^k$, 
$\dbar T_3=0$ and 
\begin{eqnarray*}
\dbar T_1&=& \lbra{ {\overline T}_1}{T_1}^{1,0}\wedge \oom^1
+\lbra{ {\overline T}_2}{T_1}^{1,0}\wedge \oom^2=\frac12 T_2\wedge\oom^1,\\
\dbar T_2 &=&\lbra{ {\overline T}_1}{T_2}^{1,0}\wedge \oom^1
+\lbra{ {\overline T}_2}{T_2}^{1,0}\wedge \oom^2=T_3\wedge\oom^1.
\end{eqnarray*}
Therefore, $\Lambda=T_2\wedge T_3$ is a holomorphic Poisson structure. By structure equations 
(\ref{poisson example structure}), $T_2\wedge T_3$ commutes with every element in $\lieg^{*(0,1)}$. 
By Corollary \ref{thm: trivial}, $DGA(M, J)=DGA(M, \cj_\Lambda)$.

\subsection{Nilmanifolds}
Suppose that $M$ is a nilmanifold $M=\Delta\backslash G$ and the generalized complex structure 
$\cj$ is invariant, there is an inclusion of left-invariant sections in the space of smooth sections 
for various bundles. 

Since the evaluation of invariant forms on invariant vectors are constants, the restriction of the 
Courant bracket to $\lieg\oplus\lieg^*$ is reduced to 
\[
\lbra{X+\alpha}{Y+\beta}=\lbra{X}{Y}+\call_X\beta-\call_Y\alpha
=\lbra{X}{Y}+\iota_X d\beta-\iota_Yd\alpha.
\]

If $\ell$ and $\oell$ represent the space of invariant sections for the Lie bialgebroid
$L$ and $\oL$ respectively, we have $\ell\oplus \oell=(\lieg\oplus\lieg^*)_\CC$ and the inclusions
\begin{equation}
\iota: \wedge^k\ell \longrightarrow C^\infty(M, \wedge^kL)
\qquad \mbox{ and } \qquad
\iota: \wedge^k\oell \longrightarrow C^\infty(M, \wedge^k\oL).
\end{equation}
Since the differential $\ddel$ on $C^\infty(M, \wedge^kL)$ is due to the invariant Lie algebroid structure
on $C^\infty(M, \oL)$, the inclusion map $\iota$ of $\ell$ intertwines with the differential $\ddel$ on $\ell$ so 
that we have an inclusion of invariant differential Gerstenhaber algebra:
\begin{equation}\label{inclusion of DGA}
\iota: DGA(\ell) \hookrightarrow DGA(M, \cj).
\end{equation}

\begin{problem}\label{p: quasi-iso} When will the inclusion map be a quasi-isomorphism?
\end{problem}

In view of Proposition \ref{proposition: symp iso} and Nomizu's Theorem \cite{Nomizu}, when the
generalized complex structure is a symplectic structure, we have an isomorphism. 
As noted in previous sections, the work regarding generalizing Nomizu's work to Dolbeault cohomology on 
complex structure on nilmanifolds remains an on-going effort in the past twenty years. 

\begin{theorem}{\rm\cite[Theorem 1]{Chen-Fino-P}} \label{t: generalized}
Suppose that $M=\Delta\backslash G$ is a nilmanifold with an abelian complex structure.
The inclusion map $\iota$ in {\rm(\ref{inclusion of DGA})} is a quasi-isomorphism.
\end{theorem}

The above theorem enables a  collection of work on generalized deformation on holomorphic Poisson 
cohomology on 
nilmanifolds \cite{Poon-Simanyi-2017} \cite{Poon-Simanyi-2019}.
In \cite{Angella-C-Kasuya} Angella et al. study
symplectic bundles over complex manifolds as generalized complex structures.
They are fiber bundles over  generalized complex manifolds such  that
a global 2-form on the total space of the bundle is restricted to a symplectic form on each fiber. 
Taking advantage of Nomizu's work and 
the quasi-isomorphism when a complex structure $J$ is abelian, 
we paraphrase a result of  \cite{Angella-C-Kasuya}. 

\begin{theorem}{\rm\cite[Theorem 5.3]{Angella-C-Kasuya}} When a generalized complex
structure on a nilmanifold $M=\Delta\backslash G$ is realized as the symplectic bundle over an abelian 
complex nilmanifold, then the inclusion map $\iota$ in {\rm(\ref{inclusion of DGA})} is a quasi-isomorphism.
\end{theorem}

\section{Semi-Abelian Generalized Complex Structure} \label{sec: semi}

We are interested in exploring ways to 
extend the concept of abelian complex structure to generalized complex manifolds. 
Consider the definition of abelian complex structure as given in (\ref{classical abelian}). 
Let $\cj$ be a generalized complex structure. One may  attempt to expand it 
na\"ively so that \it for all \rm $X, Y\in \lieg$ and $\alpha, \beta\in \lieg^*$, 
\begin{equation}\label{eq: cj abelian}
\lbra{\cj(X+\alpha)}{\cj(Y+\beta)}=\lbra{X+\alpha}{Y+\beta}.
\end{equation}
However, when $\cj$ is a symplectic structure $\Omega$ it forces 
the algebra $\lieg$ to be abelian because we would have
$$
\lbra{X}{Y}=\lbra{\Omega (X)}{\Omega (Y)}.
$$
Yet the Courant bracket between a pair of 1-forms is equal to zero. 

As our goal is to create a theory to adopt to cohomological computation to facilitate investigation deformation 
of generalized complex structures, we focus on the structure of the differential 
Gerstenhaber algebra for a generalized complex structure as given in (\ref{symbolic DGA}). 
Let us once again recall the key features of a classical complex structure $J$ being abelian is presented in the 
two equivalent conditions in Proposition \ref{prop: motivation}. This proposition and our need in computing cohomology
effectively at least in some situation drives  our proposed concept below. 

Let $\cg$ be the vector space $\lieg\oplus\lieg^*$ equipped with the 
Courant bracket induced by the Lie bracket on $\lieg$. It is a Lie algebra with the properties that $\lieg$ is a 
subalgebra and $\lieg^*$ is an abelian ideal. Therefore, $\cg$ is a semi-direct product $\cg=\lieg\rtimes \lieg^*$. 

\begin{definition}\label{def: admissible} Suppose that $\ca$ is a subalgebra of $\cg=\lieg\rtimes \lieg^*$, 
$\ck$ is an abelian ideal of 
$\cg$, and both $\ca$ and $\ck$ are maximally isotropic with respect to the natural pairing {\rm {(\ref{pairing})}}, 
the semi-direct product presentation $\cg=\ca\rtimes \ck$ is called admissible. 
We will also call $(\ca, \ck)$ an admissible pair associated to the Lie algebra $\lieg$. 
\end{definition} 

By virtual of being maximally isotropic,  the non-degenerate pairing defines an isomorphism $\ck^*\cong \ca$.

 \begin{definition}\label{def: key} Let $\lieg$ be a Lie algebra with an integrable generalized complex structure
 $\cj$.  Suppose that there exists an admissible pair $(\ca, \ck)$
 associated to $\lieg$ such that (a) both $\ca$ and $\ck$ are
 invariant of the map $\cj$ and (b) the restriction of $\cj$ on $\ca$ satisfies Identity 
 {\rm(\ref{classical abelian})}, then
 the generalized complex structure $\cj$ is said to be semi-abelian and the pair 
 $(\ca, \ck)$ is said to be $\cj$-admissible. 
 \end{definition}
 
 \begin{definition}\label{def: key on M}
Let $M=\Delta\backslash G$ be a nilmanifold with an invariant generalized complex structure $\cj$. 
It is semi-abelian if $\cj$ is semi-abelian on the Lie algebra $\lieg$ of the group $G$. 
\end{definition}

 Since $\ca$ and $\ck$ are $\cj$-invariant, their complexification have respective decomposition by 
 eigenspaces. Let $\liea$ and $\liek$ be their respective $(+i)$-eigenspaces. Their
 conjugations $\overline\liea$ and $\overline\liek$ are the $(-i)$-eigenspace so that 
 \begin{equation}
 \ell=\liea\oplus \liek, \qquad \oell=\overline\liea\oplus\overline\liek, \qquad
 \ell\oplus \oell=(\lieg\oplus\lieg^*)_\CC.
 \end{equation}

Since the restriction of $\cj$ on $\ca$ satisfies Identity (\ref{classical abelian}), $\liea$ is an abelian complex
algebra. Therefore, 
$\ell$ is a semi-direct product $\liea\rtimes \liek$ of two abelian subalgebras. As $\ca\cong \ck^*$ by the natural 
non-degenerate pairing, 
\begin{equation}\label{eq: duals}
\liea = {\overline\liek}^*, \qquad \liek= {\overline\liea}^*.
\end{equation}
Since for any element $\ell_1\in\ell$, $\ddel \ell_1$ is obtained by the evaluation on every pair of elements
$\oell_2, \oell_3$ in $ \oell$: 
\[
(\ddel \ell_1)(\oell_2, \oell_3)=-\langle \ell_1, \lbra{\oell_2}{\oell_3}\rangle,
\]
Given the structure of $\oell={\overline\liea}\rtimes {\overline\liek}$ and $\liea$ and $\liek$ are abelian, 
$\lbra{\oell}{\oell}\subseteq \overline\liek$.  By (\ref{eq: duals}), we find that $\liek$ is in the kernel of 
$\ddel$ and the image of $\liea$ via $\ddel$ is contained in $\liea\otimes \liek$. 

\begin{proposition}\label{prop: properties}
Let $\cj$ be a semi-abelian complex structure on a Lie algebra $\lieg$ with admissible pair $(\ca, \ck)$, the following
holds. 
\begin{itemize} 
\item $\liea$ is abelian subalgebra and $\liek$ is an abelian ideal such that $\ell=\liea\rtimes \liek$. 
\item $\liea = {\overline\liek}^*$ and $\liek= {\overline\liea}^*$. 
\item $\ddel\liek=\{ 0\}$ and $\ddel\liea\subseteq \liea\otimes \liek$. 
\end{itemize} 
\end{proposition}

The last point makes obvious restriction on the $DGA(\lieg, \cj)$ in terms of a lower bound on 
the dimension of its first cohomology $H^1(\lieg, \cj)$. This proposition mirrors the observation on 
abelian complex structures as noted in Proposition \ref{prop: motivation}. 
Next, we  find constraints on symplectic structure
being semi-abelian. 

\begin{proposition}\label{prop: sym} Suppose that $\Omega$ is a symplectic structure on a nilpotent algebra $\lieg$. 
If $(\ca, \ck)$ is $\Omega$-admissible, 
there exist an abelian subalgebra $\lieb$ in $\lieg$ and an abelian ideal $\lieh$ in $\lieg$ such that 
\[
\ca=\{X, \Omega(X): X\in \lieb\}, \qquad \ck=\{Y, \Omega(Y): Y\in \lieh\}
\]
with  $d\Omega(X)=0$ for all $X\in \lieh$.  In particular, $\lieg$ is a semi-direct product $\lieg=\lieb\rtimes \lieh$. 
\end{proposition}
\bproof
When $\Omega$ is a semi-abelian symplectic form, suppose that $\ell=\liea\rtimes \liek$. Then
\[
\liea=\{X-i\Omega(X): X\in \lieb\}, \qquad \liek=\{Y-i\Omega(Y): Y\in \lieh\}
\]
for some vector subspaces $\lieb$ and $\lieh$ in $\lieg$.  For any $X,Y\in \lieg$, 
\[
\lbra{X-i\Omega(X)}{Y-i\Omega(Y)}=\lbra{X}{Y}-i\Omega(\lbra{X}{Y}).
\]
Therefore, it is equal to zero if and only $\lbra{X}{Y}=0$. Therefore $\lieb$ is an abelian 
subalgebra and $\lieh$ is an abelian ideal. 
Furthermore, $\ddel (X-i\Omega(X))=0$ if and only $d\Omega(X)=0$. The dimension restriction forces $\lieg$ to be a
direct sum of $\lieb$ and $\lieh$. The Lie algebra structure follows. 
\eproof

We may now restrict the scope of Problem \ref{p: quasi-iso} to our current perspective. 

\begin{problem} Suppose that $M=\Delta\backslash G$ is a nilmanifold with a semi-abelian 
generalized complex structure. Is the inclusion map {\rm(\ref{inclusion of DGA})} necessarily 
a quasi-isomorphism? 
\end{problem}

Corollary \ref{thm: trivial} indicates that the 
concept of \it semi-abelian \rm is invariant of  deformations generated by center of 
the underlying Schouten algebra
of $DGA(M, \cj)$. 
In view of the work in \cite{CFP} \cite{MPPS-2-step} and the results as noted in Proposition \ref{prop: 6d abelian}, 
there are non-trivial deformation keeping the property of \it abelian \rm stable.  

\begin{problem} Suppose that $\cj$ is a semi-abelian generalized complex structure. Find all $\Gamma$ in 
$\wedge^2\ell$ such that $(L_{\OG}, \oL_{\Gamma})$ is a semi-abelian generalized complex structure.
\end{problem} 

\begin{remark}
The concept of semi-abelian generalized complex structure as given in Definition {\rm{\ref{def: key}}}
 on algebra
level potentially could be extended in a context similar to those in {\rm \cite{ABD} \cite{Barberis} \cite{LUV-2019}}. 
\end{remark}

\section{Examples}\label{sec: example}
By construction, all abelian complex structures are 
semi-abelian (generalized) complex structures. 
In this section, we present a collection of examples of semi-abelian generalized complex structures, 
including symplectic structures,  non-abelian nilpotent complex structures, together with
non-complex and non-symplectic
type structures. 

If a symplectic structure is semi-abelian as a generalized complex structure, we will address it as
a semi-abelian symplectic structure. 
In Example \ref{ex: 0, 0, 0, 0, 12, 13}, we find a non-abelian complex structure
such that as a generalized complex structure, it is semi-abelian. In short,  
there are semi-abelian complex structures that fail to be abelian. 
We start on four-dimension cases. 

\begin{example}\label{ex: (0,0,0,12)} The algebra $(0, 0, 0, 12)$, revisit. 
\end{example} 
Since $de^4=e^{12}$, the structure equations of Courant bracket on $\cg=\lieg\rtimes \lieg^*$ are
\begin{equation}
\lbra{e_1}{e_2}=-e_4, \qquad \lbra{e_1}{e^4}=e^2, \qquad \lbra{e_2}{e^4}=-e^1.
\end{equation}
When we choose $Je_1=e_2$ and $Je_3=e_4$, we have the well-known complex 
structure of Kodaira surface as noted in Section \ref{sec: Kodaira surfaces}. Its 
corresponding choice of  admissible pair is $\ca=\lieg$ and $\ck=\lieg^*$. 

On the other hand, we may choose 
$
\ca=\{e_1, e_2, e_4, e^3\}$  and 
$\ck=\{e^1, e^2, e^4, e_3\}.$
Consider a type-1 generalized complex structure $\cj$ defined by its components.  
\[
Je_1=e_2, \quad Je_2=-e_1, \quad B=e^{34}, \quad \Pi=e_{34}. 
\]
In particular, $\cj(e_3)=e^4, \cj(e_4)=-e^3$. 
The pair $(\ca, \ck)$ becomes $\cj$-admissible. Moreover, the $(+i)$-eigenspace is spanned by 
the following elements. 
\[
\ell_1=e_1-ie_2, \quad \ell_2=e_4+ie^3, \quad \ell_3=e^1-ie^2, \quad \ell_4=e_3-i e^4.
\]
Since $e_4, e^1, e^2, e^3$ are in the center of the algebra $\cg$, the integrability is due to 
\[
 \lbra{\ell_1}{\ell_4}=\lbra{e_1-ie_2}{e_3-i e^4} =-i\lbra{e_1-ie_2}{ e^4}  = -i(e^2+ie^1)=\ell_3. 
 \]
Therefore, $\cj$ is integrable  and  semi-abelian. 

 \begin{example}\label{ex: counter} An example admitting no semi-abelian generalized complex structures.
 \end{example} 
Consider the algebra $(0, 0, 12, 13)$. 
Since $de^3=e^{12}$, $de^4=e^{13}$, it admits  symplectic structures such as $e^{23}+e^{14}$. 
However, this algebra does allow a semi-direct product as prescribed by Proposition \ref{prop: sym}, 
there is no semi-abelian symplectic structure on this algebra. 

The structure of Courant bracket on $\cg=\lieg\rtimes \lieg^*$ are 
\begin{eqnarray*}
& \lbra{e_1}{e_2}=-e_3, \qquad \lbra{e_1}{e_3}=-e_4 &\\
&\lbra{e_1}{e^3}=e^2, \qquad \lbra{e_1}{e^4}=e^3, \qquad
\lbra{e_2}{e^3}=-e^1, \qquad \lbra{e_3}{e^4}=-e^1. &
\end{eqnarray*}
The center of $\cg$ is spanned by $e_4, e^1, e^2$. 

Utilizing  
ascending basis similar to (\ref{generic constraint}) as in \cite{Malcev} \cite{Salamon}, Cavalcanti et al.
provided a characterization of generalized complex structures in terms of spinor formalism on nilmanifolds
 \cite[Corollary 1]{Cavalcanti-G}. 
In particular, when $\cj$ is type-1, there exists a spinor representation 
$\rho=e^{B+i\om}\theta$ with real 2-forms $B$ and $\om$, and 
complex 1-form $\theta$ such that $\om\wedge \theta\wedge \otheta\neq 0$ and $d\theta=0$. 
Given the algebra at hand, 
the only choice for $\theta$, up to linear combination, is $\theta=e^1+ie^2$. 
It follows that up to a constant 
\[
\om=e^{34}+e^3\wedge (a_{31}e^1+a_{32}e^2)+e^4\wedge (a_{41}e^1+a_{42}e^2)
\]
for some real numbers $a_{31}, a_{32}, a_{41}, a_{42}$. 
Therefore the generalized complex structure $\cj$ is given by
\begin{eqnarray*}
&\cj e_1=e_2,  \quad  \quad \cj e^1=e^2, &\\ 
& \cj e_3=e^4+a_{31}e^1+a_{32}e^2, \qquad \cj e_4=-e^3+a_{41}e^1+a_{42}e^2&.
\end{eqnarray*}
It follows that the $(+i)$-eigenspace is spanned as 
\[
\ell=\{ e_1-ie_2, e^1-ie^2, e_3-i(e^4+a_{31}e^1+a_{32}e^2), e_4-i(-e^3+a_{41}e^1+a_{42}e^2) \}.
\]
As $\ell$ is isotropic, 
$a_{31}-ia_{32}=0$, $ a_{41}-ia_{42}=0.$
Therefore, all such constants are equal to zero so that
\begin{equation}
\ell= \{ e_1-ie_2, e^1-ie^2, e_3-ie^4, e_4+ie^3 \}.
\end{equation}
As $e^1$ and $e^2$ are in the center of the algebra $\cg$, to verify integrability of $\cj$ it suffices to see that 
\begin{eqnarray*}
&& \lbra{e_1-ie_2}{e_3-ie^4}=\lbra{e_1}{e_3-ie^4}-i\lbra{e_2}{e_3-ie^4}= -e_4-ie^3=-(e_4+ie^3), \\
&&\lbra{e_1-ie_2}{e_4+ie^3}=i\lbra{e_1}{e^3}+\lbra{e_2}{e^3}=ie^2-e^1=-(e^1-ie^2).
\end{eqnarray*}
In terms of $\oell$,
$\lbra{e_1+ie_2}{e_3+ie^4}=-(e_4-ie^3)$, 
$\lbra{e_1+ie_2}{e_4-ie^3}=-(e^1+ie^2).$
 It follows that no combination of $e_3-ie^4$ and $e_1-ie_2$ is in the kernel of the operator $\ddel$. 
Moreover, 
\[
\ker\ddel =\{ e^1-ie^2, e_4+ie^3\}.
\]
Now suppose that $(\ca, \cg)$ is a $\cj$-admissible pair. By Proposition \ref{prop: properties},
 $\liek$ is in the kernel of $\ddel$. 
A dimension restriction shows that  $\ck=\{e_4, e^1, e^2, e^3\}$. 
While $\ck$ is indeed an abelian ideal in $\cg$, the dual space is
$
\ck^*=\ca=\{e^4, e_1, e_2, e_3\}, 
$
and it fails to be a Lie subalgebra. Therefore, the generalized complex structure $\cj$ is not semi-abelian. 

Finally by \cite[Theorem 3.2]{Cavalcanti-G} as well as \cite[Proposition 2.1]{Snow}, 
the algebra $(0,0,12, 13)$ does not admit any invariant complex structure, let alone a semi-abelian one. 

It is an elementary fact that  the only non-abelian 4-dimensional nilpotent algebras 
are $(0, 0, 0, 12)$ and $(0, 0, 12, 13)$, see e.g., \cite{Magnin}. The observation above shows that 
the latter admits type-1 as well as type-0 generalized complex structure, but none of them are 
semi-abelian. 

\begin{proposition} The only  four-dimensional  nilpotent algebras admitting 
semi-abelian generalized complex structures are $(0, 0, 0, 0)$ and  $(0, 0, 0, 12)$. 
\end{proposition}

\begin{example}\label{ex: 12+34}  A type-2 example on  $(0,0,0,0,0, 12+34)$.
\end{example}
As the structure for Lie algebra $\lieg$ is  $de^6=e^{12}+e^{34}$, 
the structural equations associated to $\cg=\lieg\rtimes\lieg^*$ are
\begin{eqnarray*}
&\lbra{e_1}{e_2}=-e_6, \qquad \lbra{e_3}{e_4}=-e_6,& \\
&\lbra{e_1}{e^6}=e^2, \quad  \lbra{e_2}{e^6}=-e^1, \quad \lbra{e_3}{e^6}=e^4, \quad \lbra{e_4}{e^6}=e^3.&
\end{eqnarray*} 
As seen in Example \ref{ex: Kodaira},  $(\lieg, \lieg^*)$ forms an admissible pair for 
an abelian complex structure. 
On the other hand, the pair below is also admissible.
\[
\ca = \{e_1, e_2, e_3, e_4, e_6, e^5\}, \quad 
\ck=\{e^1, e^2, e^3, e^4, e^6, e_5\}. 
\]
Consider a generalized complex structure $\cj$ given below. 
\begin{eqnarray*}
&\cj e_1=e_2, \cj e_2=-e_1, \cj e_3=e_4, \cj e_4=-e_3,  \cj e_6=-e^5, \cj e^5=e_6, & \\
&\cj e^1=e^2, \cj e^2=-e^1, \cj e^3=e^4, \cj e^4=-e^3, \cj e_5=e^6,\cj e^6=-e_5.
\end{eqnarray*}
It leaves $\ca$ and $\ck$ invariant. 
In terms of components in the matrix representation for $\cj$,
\begin{equation}
J=e^1\otimes e_2-e^2\otimes e_1+e^3\otimes e_4-e^4\otimes e_3, 
\quad B=e^5\wedge e^6, \quad \Pi=e_5\wedge e_6.
\end{equation}
The $(+i)$-eigenspace for $\cj$ is
\[
\ell=\{e_1-ie_2, e_3-ie_4, e^1-ie^2, e^3-ie^4, e_5-ie^6, e_6+ie^5\}.
\]
Since $de^1=de^2=de^3=de^4=0$, $e^1-ie^2, e^3-ie^4$ are in the center of the Schouten bracket on $\ell$. 
Similarly, $e_6$ is in the center of $\lieg$ and $de^5=0$, $e_6+ie^5$ is also in the center of the Schouten bracket 
on $\ell$. 
Therefore, a computation of Schouten bracket on $\ell$ is reduced to 
\begin{eqnarray*}
&&\lbra{e_1-ie_2}{e_3-ie_4}=0,\\
&&\lbra{e_1-ie_2}{e_5-ie^6}=-i\lbra{e_1-ie_2}{e^6}=-i(e^2+ie^1)=e^1-ie^2.\\
&&\lbra{e_3-ie_4}{e_5-ie^6}=i\lbra{e_3-ie_4}{e^6}=-i(e^4+ie^3)=e^3-ie^4.
\end{eqnarray*}
Since $\ell$ is closed with respect 
to the Schouten bracket, $\cj$ is an integrable generalized complex structure. Moreover, 
$\ell=\liea\rtimes \liek$ where 
\begin{equation}
\liea=\{e_1-ie_2, e_3-ie_4,  e_6+ie^5\}, 
\quad
\liek=\{ e^1-ie^2, e^3-ie^4, e_5-ie^6\}.
\end{equation}
It follows that $\cj$ is a type-2 semi-abelian complex structure. 
Since $de^5=0$ and $de^6=e^{12}+e^{34}$, they are type-(1,1) when treated as
2-form on the quotient space spanned by $\{e_1, e_2, e_3, e_4\}$, 
we obtain an example for Proposition \ref{prop: toric} when we consider 
taking the quotient of the algebra $(0,0, 0,0, 0, 12+34)$ by its center as a principal bundle
projection on the nilmanifold level. 

 \begin{example}\label{ex: 0, 0, 0, 0, 12, 14+25}  A semi-abelian symplectic structure on $(0, 0, 0, 0, 12, 14+25)$.
 \end{example}
Let $\lieb=\{e_2, e_3, e_4\}$ and $\lieh=\{e_1, e_5, e_6\}$, then $\lieg=\lieb\rtimes \lieh$.
This algebra   has a symplectic form, 
$
\Omega=e^{13}+e^{26}+e^{45}.$  Since 
$\Om (e_1)=e^3$, $\Om (e_5)=-e^4$, and $\Om (e_6)=-e^2$, the pair $(\lieb, \lieh)$ satisfies 
 the 
conditions in Proposition \ref{prop: sym}  with respect to $\Omega$.
Therefore $\Omega$ is an abelian symplectic structure with $\ell=\liea\rtimes \liek$ where
\[
\liea=\{e_2-ie^6,  e_3+ie^1,  e_4-ie^5 \}, \qquad 
\liek=\{e_1-ie^3,   e_5+ie^4,  e_6+ie^2\}.
\]

\begin{example}\label{ex: 0, 0, 0, 0, 12, 13} A semi-abelian, non-abelian complex structure on $(0, 0, 0, 0, 12, 13)$.
\end{example} 
The structure
equations for  this $\cg=\lieg\rtimes \lieg^*$ are given below. 
\begin{eqnarray*}
&\lbra{ e_1}{e_2}=-e_5, \qquad \lbra{e_1}{e_3}=-e_6, & \\
& \lbra{e_1}{e^5}=e^2, \quad \lbra{e_2}{e^5}=-e^1, \quad \lbra{e_1}{e^6}=e^3, \quad \lbra{e_3}{e^6}=-e^1. &
\end{eqnarray*}
Although the algebra $\lieg$  does not admit
any abelian complex structure \cite[Proposition 3.3]{CFU},
it does admit nilpotent complex structures \cite{Salamon}. We adopt the one given 
in \cite[Table 1]{Cavalcanti-G}. Namely, $\cj=J$, $B=0$, $\Pi=0$, with
\begin{equation}\label{nilpotent eg}
J e_1=e_4, \quad J e_2=e_3, \quad J e_5=e_6.
\end{equation} 
The pair
$\ca=\{e_2, e_3, e^1, e^4, e^5, e^6\}$, 
$\ck=\{e_1, e_4, e_5, e_6, e^2, e^3\}$
 is $\cj$-admissible. Therefore, the nilpotent complex structure in (\ref{nilpotent eg}) 
is semi-abelian although it is not abelian. 

\begin{example}\label{ex: 0, 0, 0, 0, 12, 13}
A semi-abelian holomorphic Poisson structure on 
$(0,0,0,0, 12, 13)$.
\end{example}
A basis for $\lieg^{1,0}$ with respect to the complex structure (\ref{nilpotent eg})  is 
\begin{equation}
T_1=\frac12(e_1-ie_4), \qquad T_2=\frac12(e_2-ie_3), \qquad T_3=\frac12(e_5-ie_6).
\end{equation}
It follows that  
$
\oom^1=e^1-ie^4$,  $ \oom^2=e^2-ie^3$,  and $\oom^3=e^5-ie^6$
forms a basis for 
$\lieg^{*(0,1)}$. Moreover, $d\oom^1=0$,  $d\oom^2=0$, and
\[
d\oom^3=\frac12 \om^1\wedge\oom^2+\frac12\oom^1\wedge\oom^2. 
\]
Equivalently, 
\begin{equation}
\lbra{T_1}{T_2}= -\frac12 T_3, \quad \lbra{{\overline T}_1}{T_2} =-\frac12 T_3.
\end{equation}
By Proposition \ref{prop: classical},   $\dbar T_2=-\frac12 T_3\wedge \oom^1$. 
It follows that $\Lambda=T_2\wedge T_3$ is a holomorphic Poisson bivector. 
 Moreover, $\lbra{T_2}{\oom^3}=0$ and $\lbra{T_3}{\oom^3}=0$,  
  $ad_\Lambda=0$.  By Corollary \ref{thm: trivial},
 $\cj_\Lambda$ defines a semi-abelian generalized complex structure.

\

\noindent{\textbf{Acknowledgement.}} The author thanks 
 Zhuo Chen and Anna Fino for their helpful suggestions
on the initial draft of this manuscript. He thanks the referee for helpful suggestions. 
He is very grateful for Sau King Chiu's support during   
the pandemic when this piece of work was conceived.

\end{document}